\documentclass[11pt,a4paper]{article}
\usepackage[T1]{fontenc}
\usepackage[utf8]{inputenc}
\usepackage[english]{babel}
\usepackage{lmodern,amsmath,amssymb,amsthm,mathtools}
\usepackage[a4paper,margin=27mm]{geometry}
\usepackage{microtype,booktabs,tabularx}
\usepackage{enumitem}
\usepackage{tikz}
\usetikzlibrary{arrows.meta,positioning,fit}
\usepackage[hidelinks]{hyperref}
\usepackage{fancyhdr}
\newtheorem{theorem}{Theorem}
\newtheorem{lemma}[theorem]{Lemma}
\newtheorem{proposition}[theorem]{Proposition}
\newtheorem{corollary}[theorem]{Corollary}
\theoremstyle{definition}
\newtheorem{remark}[theorem]{Remark}
\newcommand{\reg}{\operatorname{reg}}
\newcommand{\Ind}{\operatorname{Ind}}
\newcommand{\sd}{\operatorname{sd}}
\newcommand{\lk}{\operatorname{lk}}
\newcommand{\del}{\operatorname{del}}
\newcommand{\Ht}{\widetilde H}
\newcommand{\ZZ}{\mathbb Z}
\newcommand{\NN}{\mathbb Z_{\geq 0}}
\newcommand{\m}{\mathfrak m}
\newcommand{\G}{\Gamma}
\title{\bfseries A proof of the Nevo--Peeva conjecture: Exact regularity and linear powers of gap-free edge ideals}
\author{Luis Alfredo Dupont Garc\'ia\\[3pt]
\small Faculty of Mathematics, Universidad Veracruzana\\
\small Xalapa, Veracruz, Mexico\\[3pt]
\small E-mail: \href{mailto:ldupont@uv.mx}{ldupont@uv.mx}\\
\small ORCID: \href{https://orcid.org/0000-0002-9044-3547}{0000-0002-9044-3547}}
\date{}
\hypersetup{pdfauthor={Luis Alfredo Dupont Garc\'ia},pdftitle={A proof of the Nevo--Peeva conjecture: Exact regularity and linear powers of gap-free edge ideals},pdfsubject={Commutative algebra and combinatorics},pdfkeywords={edge ideals, gap-free graphs, regularity, Nevo-Peeva conjecture, linear powers}}
\begin{document}
\maketitle
\thispagestyle{plain}
\begin{abstract}
Let $G$ be a finite simple graph with at least one edge, and let $I(G)$ be its edge ideal over a field $K$. We prove that if $G$ is gap-free, then
\[
 \reg_K I(G)^q=\max\{2q,\reg_K I(G)+q-1\}\qquad(q\geq1).
\]
This answers both open problems of Nevo and Peeva on ordinary powers: $I(G)^q$ has a linear resolution exactly when $q\geq\reg_K I(G)-1$, and their proposed regularity recurrence holds with equality. As consequences, the upper bound $\reg I(G)^q\leq 2q+\reg I(G)-2$ conjectured by Alilooee, Banerjee, Beyarslan, H\`a, Jayanthan and Selvaraja holds for gap-free graphs; every power $I(G)^q$ with $q\geq2$ has a linear resolution when $\reg_K I(G)=3$; and the stabilization index equals $\reg_K I(G)-1$, which is unbounded on this class. The proof rests on a combinatorial result of independent interest: a parity-sensitive homology vanishing theorem for a filtration of the independence complex by a weighted Hall defect. An internalization lemma and sublevel complexes of submodular functions control the local terms, and Hall's theorem together with an implication-digraph argument classifies the zero level. A capacitated form of the Tutte--Berge formula and Takayama's formula identify the filtration with the degree complexes of ordinary powers. The algebraic statements hold over every field, and the topological vanishing holds with arbitrary abelian coefficients.
\end{abstract}

\section{Introduction}
Let $G=(V,E)$ be a finite simple graph, let $S=K[x_v:v\in V]$ be a standard graded polynomial ring over a field $K$, and write
\[
 I=I(G)=(x_ux_v:uv\in E)\subseteq S.
\]
Throughout the main algebraic statements we assume $E\ne\varnothing$. Regularity is taken for the ideal as a graded $S$-module, so that $\reg I=\reg(S/I)+1$. A graph is \emph{gap-free} if it has no induced subgraph consisting of two disjoint edges. Equivalently, its complement has no induced cycle of length four. That gap-freeness is necessary for some ordinary power to have a linear resolution is due to Francisco, H\`a and Van Tuyl (recorded in \cite[Proposition 1.8]{NP}; see also Section~\ref{sec:consequences}). We use the spelling ``gap-free'' throughout; ``gapfree'' is also common.

Fr\"oberg's theorem \cite{Fro} characterizes the case $\reg I=2$ by chordality of the complement; Herzog, Hibi and Zheng \cite{HHZ} showed that in this case all powers have linear resolutions. For a survey of regularity bounds and exact formulas, see \cite{BBHsurvey}. The asymptotic linearity theorems of Cutkosky--Herzog--Trung and Kodiyalam \cite{CHT,Kod} place the question of linear resolutions in the study of the function $q\mapsto\reg I^q$. For an edge ideal it eventually has the form $2q+b(I)$ with an integer $b(I)\geq0$.

Nevo and Peeva \cite[Open Problems 1.11(1) and (2)]{NP} asked whether gap-freeness is equivalent to eventual linearity of the powers and, more strongly, whether for gap-free graphs the inequality
\[
 \reg I^{q+1}\leq \max\{2q+2,\reg I^q+1\}\qquad(q\geq1)
\]
holds. In the discussion immediately after Open Problems 1.11, condition (b) specifies the predicted range $q\geq\reg I-1$. Subsequent work, including \cite{BBH,BN,E26}, refers to these problems as the Nevo--Peeva conjecture. Nevo and Peeva's Counterexample 1.10 shows that gap-freeness alone does not force the square to be linear. Numbering refers to the published article \cite{NP}, which differs from that of the author manuscript.

Partial results include Nevo's theorem for gap-free claw-free graphs \cite{Nevo}, Banerjee's theorem for gap-free cricket-free graphs \cite{Ban}, the local-to-global bounds of \cite{BBH}, and the extension approach of \cite{Kanno}; the induced-matching lower bound of Beyarslan--H\`a--Trung \cite{BHT} and Erey's diamond-free case \cite{Erey18} provide further context. Banerjee and Nevo \cite{BN} proved the general upper-bound conjecture for squares and for all powers of bipartite graphs; their discussion of the dunce hat identifies an obstruction to a direct colon-ideal induction. Minh and Vu \cite[Theorem 1.3]{MV} proved the exact regularity formula for the second and third powers, and their Conjecture A \cite[p.~339]{MV} asks for the equality
\[
 \reg I^q=\reg I^{(q)}=\max\{2q,\reg I+q-1\}\qquad(q\geq2)
\]
for gap-free graphs, where $I^{(q)}$ denotes the symbolic power. We address only its assertion for ordinary powers. More recent work on linear quotients includes \cite{E26,E3}; the separate preprint \cite{D4} gives a linearity criterion for ordinary powers and applies it to gap-free graphs of matching number at most four. None of these last three results is used in the proof below. Our main result is uniform in the exponent and in the initial regularity.

\begin{theorem}\label{thm:main}
Let $G$ be a finite simple gap-free graph with at least one edge, and let $K$ be any field. For every integer $q\geq1$,
\begin{equation}\label{eq:main}
 \boxed{\reg_K I(G)^q=\max\{2q,\reg_K I(G)+q-1\}.}
\end{equation}
Consequently, $I(G)^q$ has a linear free resolution if and only if $q\geq \reg_K I(G)-1$.
\end{theorem}
Theorem~\ref{thm:main} proves the ordinary-power assertion of \cite[Conjecture A]{MV}. The equality holds for each fixed field; we do not claim that the initial regularity is independent of the characteristic. Symbolic powers and linear quotients are not addressed.

\subsection{The combinatorial results and the proof strategy}\label{subsec:strategy}
We state the two other principal results in advance. Let $a:V\to\NN$ and $A=\sum_v a_v$. Write $\Ind(G)$ for the independence complex. For an independent set $Z$ define its \emph{Hall defect} $L(Z)=a(N_G(Z))-a(Z)$, and for an independent $F$ define
\[
 \lambda(F)=\min_{Z\in\Ind(G),\,Z\supseteq F}L(Z),\qquad
 \Gamma_c=\{F\in\Ind(G):\lambda(F)\leq c\}.
\]
We call the increasing family $(\Gamma_c)_{c\in\ZZ}$ of subcomplexes of $\Ind(G)$ the \emph{Hall filtration}; the name refers to the neighborhood inequalities of Hall's theorem, which govern its zero level (Lemma~\ref{lem:zero}). Put $\epsilon(A,c)=0$ if $A+c$ is odd and $1$ otherwise.

\medskip\noindent\textbf{Combinatorial vanishing (Theorem~\ref{thm:C}).}
\emph{If $G$ is gap-free, $c\in\ZZ$ and $A\geq c+1$, then $\Ht_j(\Gamma_c;B)=0$ for all $j\geq c+\epsilon(A,c)$ and every abelian group $B$.}

\medskip\noindent\textbf{Degree-complex identity (Proposition~\ref{prop:degree}).}
\emph{If $G$ is gap-free, $q\geq1$ and $a\in\NN^V$, then the multigraded degree complex of the ordinary power satisfies}
\[
 \Delta_a(I(G)^q)=\Gamma_{2q-1-A}(G,a).
\]
Thus for fixed $a$, all levels $c\equiv A+1\pmod2$ with $c\geq1-A$ are exactly the levels occurring for positive ordinary powers. Levels of that parity below $1-A$ are void, since $L(Z)\geq-A$. The opposite parity is needed in the induction. The proposition also reduces negative multidegrees to an induced residual graph and a smaller exponent. Thus the Hall filtration is not an auxiliary device: it models all degree complexes of ordinary powers.

\paragraph{Principal consequences.}
With $R=\reg_K I(G)$, the main theorem yields the upper bound $\reg I^q\leq2q+R-2$ conjectured by Alilooee, Banerjee, Beyarslan, H\`a, Jayanthan and Selvaraja (ABBHJS), the full base case $R\leq3$ of Nevo--Peeva for $q\geq2$, the exact stabilization index $R-1$, and the criterion for the square to be linear. These consequences, the unboundedness of the stabilization index, and the fractional matching interpretation are proved in Sections~\ref{sec:algebra} and \ref{sec:consequences}.

\begin{figure}[htbp]
\centering
\begin{tikzpicture}[x=1cm,y=1cm,>=Stealth, every node/.style={font=\small},
 block/.style={draw=black!65,rounded corners=2pt,align=center,text width=3.15cm,minimum height=1.05cm,inner sep=4pt}, arr/.style={->,semithick,draw=black!75}]
 \node[block] (a) at (0,0) {Lemma~\ref{lem:bi}\\Independent-set inequality};
 \node[block] (b) at (4,0) {Lemma~\ref{lem:internal}\\Internalization};
 \node[block] (d) at (8,0) {Lemmas~\ref{lem:submod}, \ref{lem:local}\\Submodular vanishing};
 \node[block] (e) at (0,-1.85) {Lemmas~\ref{lem:negative}, \ref{lem:zero}\\Negative and zero levels};
 \node[block] (f) at (4,-1.85) {Theorem~\ref{thm:relative}\\Relative vanishing};
 \node[block] (g) at (8,-1.85) {Theorem~\ref{thm:C}\\Absolute vanishing};
 \node[block] (h) at (0,-3.70) {Lemma~\ref{lem:matching}\\Capacitated matchings};
 \node[block] (i) at (4,-3.70) {Proposition~\ref{prop:degree}\\Degree complexes};
 \node[block] (j) at (8,-3.70) {Theorem~\ref{thm:main}\\Exact regularity};
 \draw[arr] (a)--(b); \draw[arr] (b)--(d); \draw[arr] (a)--(e);
 \draw[arr] (e)--(f); \draw[arr] (d.south) -- ++(0,-.32) -| (f.north);
 \draw[arr] (f)--(g); \draw[arr] (h)--(i); \draw[arr] (i)--(j);
 \draw[arr] (g)--(j);
 \node[font=\footnotesize] at (6,-4.56) {Algebraic step uses \v{C}ech--Takayama local cohomology.};
\end{tikzpicture}
\caption{Main dependencies of the proof. The low-level results also enter the absolute induction, and the clone construction in Lemma~\ref{lem:zero} is reused in Lemma~\ref{lem:matching}.}\label{fig:dependencies}
\end{figure}
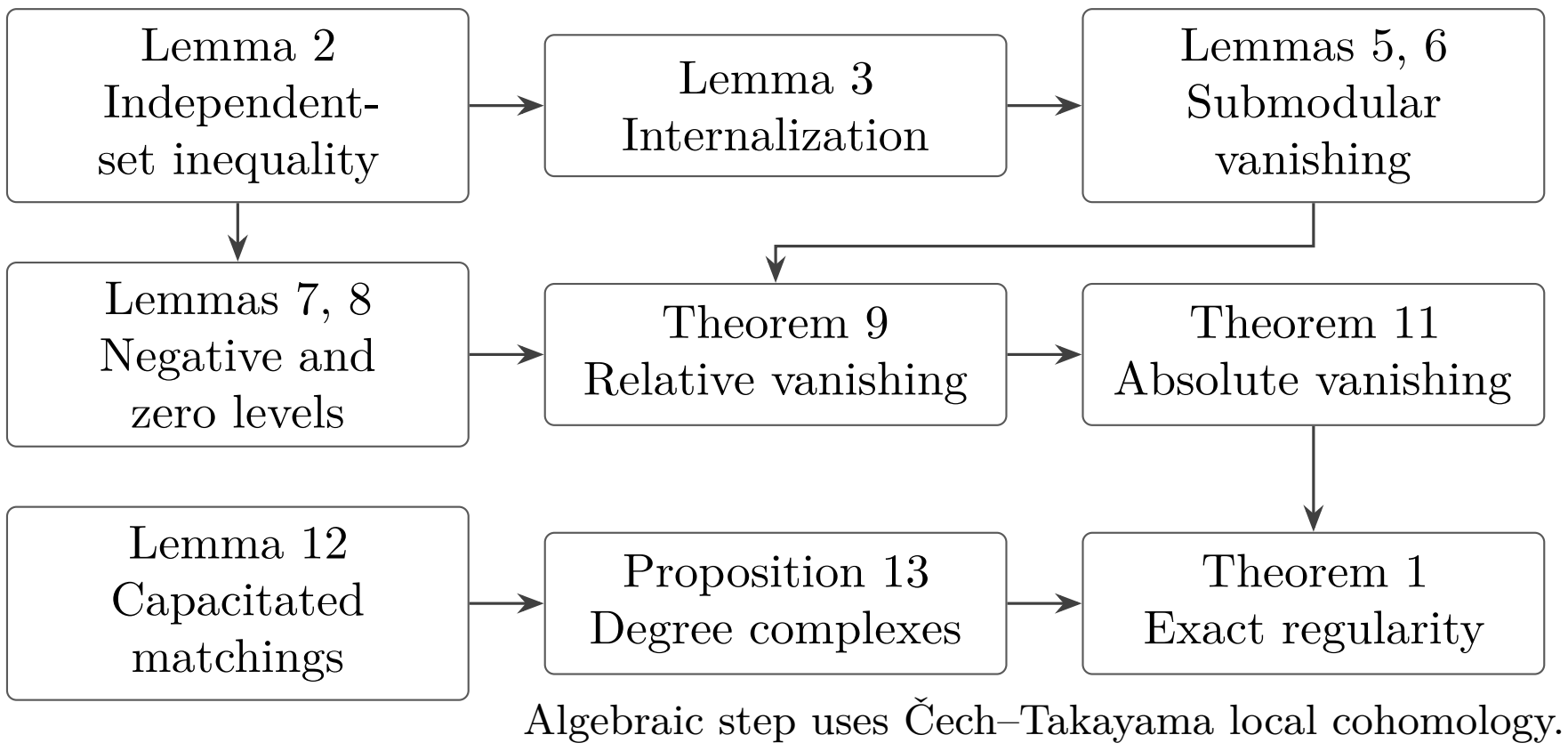

The local step is Lemma~\ref{lem:internal}: a witness at the preceding integer level can be replaced by one within an exact set. Lemma~\ref{lem:local} then follows from submodularity. The global step combines the local bound with the classification of the zero level, using relative homology of successive levels and a separate deletion/link induction for zero weights. Finally, the matching identity supplies the degree-complex formula; the parity-sensitive bound controls high-mass multidegrees in Takayama's formula. The low-mass contribution is compared with the first power. Figure~\ref{fig:dependencies} summarizes this organization.

\section{Weighted independent sets and internalization}\label{sec:prelim}
For a graph $G$ (not necessarily gap-free), let $N_G(X)$ denote the set of vertices adjacent to at least one vertex of $X$, and let $N_G[X]=X\cup N_G(X)$. We omit the subscript when the graph is clear. The graph induced by $V\setminus N_G[X]$ is denoted by $G-N_G[X]$. Let $a\in\NN^V$ be an integer weight vector and put $a(X)=\sum_{v\in X}a_v$, $A=a(V)$. For independent sets, define
\begin{align*}
 L_G(Z)&=a(N_G(Z))-a(Z),\\
 \lambda_G(F)&=\min\{L_G(Z): Z\in\Ind(G),\ F\subseteq Z\},\\
 \mu_G&=\lambda_G(\varnothing).
\end{align*}
Here $\Ind(G)$ is the independence complex, and the minimum defining $\lambda_G$ is taken only for independent $F$. For an integer $c$, the Hall filtration is defined by
\[
 \G_c(G,a)=\{F\in\Ind(G):\lambda_G(F)\leq c\}.
\]
This is a simplicial complex since $\lambda_G$ is nondecreasing under inclusion. In particular, $\mu_G\leq0$ since $L_G(\varnothing)=0$. When taking an induced subgraph we always restrict the weights to that graph; when restricting a simplicial complex to a vertex subset, its original faces and weights are retained. These two operations need not agree. The notation used throughout the proof is collected below; the least minimizer $K_0$ is shown to exist after Lemma~\ref{lem:bi}.
\begin{center}
\small
\begin{tabularx}{\linewidth}{@{}lX@{}}
\toprule
Symbol or term & Meaning\\
\midrule
$A$, $a(X)$ & Total mass; mass of a subset.\\
$L(Z)$, $\lambda(F)$ & Hall defect; minimum defect over independent extensions of $F$.\\
$\mu_G$, $K_0$ & Global minimum; least independent minimizer under inclusion.\\
$H_F$ & Induced residual $G-N_G[F]$ with restricted weights.\\
$\Gamma_c$, $\epsilon(A,c)$ & Level $c$; zero for $A+c$ odd, one for $A+c$ even.\\
Exact, tight & $L(F)=\lambda(F)=c$ (at level $c$); $L(F)=0$ (when $\mu_G=0$).\\
\bottomrule
\end{tabularx}
\end{center}

We use augmented simplicial chains with coefficients in an arbitrary abelian group $B$. The void complex (having no faces) has zero chain groups in all degrees, whereas the complex $\{\varnothing\}$ has reduced homology $B$ in degree $-1$ and zero in all other degrees. In a join, $\{\varnothing\}$ is the identity. Relative homology in this paper is the homology of the quotient of the augmented chain complexes. The notation $\G_c[X]$ denotes the induced subcomplex on $X$; $\lk_F D=\{T\subseteq V(D)\setminus F:F\cup T\in D\}$ and $\del_v D=\{T\in D:v\notin T\}$ denote link and deletion (with ambient vertices understood). For a finite complex $D$ with at least the empty face, $\sd D$ denotes its barycentric subdivision, namely the order complex of its nonempty faces; this definition retains $\{\varnothing\}$ when there are no vertices. All uses of geometric contractibility refer to complexes with vertices. The four cases used repeatedly below have the following augmented reduced homology:
\begin{center}
\begin{tabular}{@{}lcc@{}}
\toprule
Complex & $\Ht_{-1}(-;B)$ & $\Ht_0(-;B)$ \\
\midrule
Void (no faces) & $0$ & $0$ \\
$\{\varnothing\}$ & $B$ & $0$ \\
Point & $0$ & $0$ \\
$S^0$ & $0$ & $B$ \\
\bottomrule
\end{tabular}
\end{center}
All other reduced homology groups of these four complexes vanish.

\begin{lemma}[Independent-set inequality]\label{lem:bi}
For independent sets $X,Y$, define
\[
 X\wedge Y=X\cap Y,\qquad
 X\vee Y=(X\setminus N(Y))\cup(Y\setminus N(X)).
\]
Both are independent. These operations need not form a lattice on the independent sets; we use only the inequality
\begin{equation}\label{eq:bi}
 L(X\wedge Y)+L(X\vee Y)\leq L(X)+L(Y).
\end{equation}
\end{lemma}
\begin{proof}
The second operation removes both ends of every edge between $X$ and $Y$, hence yields an independent set. For an independent $Z$, a vertex contributes $-a_v$, $a_v$, or $0$ to $L(Z)$ according as it lies in $Z$, $N(Z)$, or outside $N[Z]$. Compare the normalized contributions (the signs $-1,0,1$ before multiplication by $a_v$) at a vertex $v$. If $v\in X\cap Y$, both sides of \eqref{eq:bi} contribute $-2$. If $v\in X\setminus Y$ has no neighbor in $Y$, both sides contribute $-1$. If $v\in X\setminus Y$ has a neighbor in $Y$, it is removed from $X\vee Y$; all its neighbors in $Y$ are removed as well. It has no neighbor in $X$ or in $X\cap Y$. Both sides therefore contribute zero. The case $v\in Y\setminus X$ is symmetric. Finally, for $v\notin X\cup Y$, membership in $N(X\cap Y)$ requires membership in both $N(X)$ and $N(Y)$, while membership in $N(X\vee Y)$ requires membership in at least one. Thus the left contribution never exceeds the right contribution. Multiplying by the nonnegative weights and summing proves the inequality.
\end{proof}
In particular, the intersection of all global minimizers of $L$ is itself a minimizer, because the two terms on the left of \eqref{eq:bi} are bounded below by $\mu_G$. Denote this least minimizer by $K_0$. Also, if $X,Y$ contain a fixed independent set $F$, then both $X\wedge Y$ and $X\vee Y$ contain $F$.

For an independent $F$, let $H_F=G-N[F]$, equipped with restricted weights. Every independent extension of $F$ is uniquely $F\cup T$ for $T\in\Ind(H_F)$. Since $N_G(F\cup T)=N_G(F)\sqcup N_{H_F}(T)$, we have
\begin{align}
 L_G(F\cup T)&=L_G(F)+L_{H_F}(T),\label{eq:ext}\\
 \lambda_G(F\cup T)&=L_G(F)+\lambda_{H_F}(T),\nonumber\\
 \lambda_G(F)&=L_G(F)+\mu_{H_F},\label{eq:lambdaext}\\
 \lk_F\G_c(G,a)&=\G_{c-L_G(F)}(H_F,a|_{V(H_F)}),\label{eq:link}\\
 a(V(H_F))&=A-2a(F)-L_G(F).\label{eq:mass}
\end{align}
The link identity holds also for a nonface if its link is defined by the condition $F\cup T\in\G_c$, in which case both sides are void.

\begin{lemma}[Internalization]\label{lem:internal}
Suppose that an independent set $z$ is exact at level $c$, that is, $L(z)=\lambda(z)=c$. Then
\begin{equation}\label{eq:internal}
 \G_{c-1}[z]=\langle Y\subseteq z:L(Y)\leq c-1\rangle,
\end{equation}
where the brackets mean the complex generated by the indicated sets. No gap-free hypothesis is required.
\end{lemma}
\begin{proof}
Put $H=G-N[z]$. By \eqref{eq:lambdaext}, exactness gives $\mu_H=0$. For any independent set $W$, write $J=z\setminus N(W)$ and $T=W\cap V(H)$. Then $z\vee W=J\sqcup T$. The sets $N_G(J)$ and $N_H(T)$ lie, respectively, in $N_G(z)$ and $V(H)$; they are disjoint subsets of $N_G(J\cup T)$. Consequently
\[
 L_G(z\vee W)\geq L_G(J)+L_H(T)\geq L_G(J).
\]
By Lemma~\ref{lem:bi},
\begin{equation}\label{eq:two}
 L(z\cap W)+L(J)\leq c+L(W).
\end{equation}
If $F\subseteq z$ lies in $\G_{c-1}$, choose an independent $W\supseteq F$ with $L(W)\leq c-1$. Both $z\cap W$ and $J$ contain $F$ (because $W$ is independent). Their integer defects have sum at most $2c-1$. At least one thus has defect at most $c-1$, and is a witness contained in $z$. This proves one inclusion in \eqref{eq:internal}; the reverse follows immediately from the definition of $\lambda$.
\end{proof}
\begin{remark}\label{rem:limits}
The preceding lemma concerns exact sets and the immediately preceding integer level. Neither condition may be dropped in general; see Section~\ref{subsec:sharpness}.
\end{remark}

\section{Local homology of submodular sublevels}\label{sec:submodular}
A set function $g$ is submodular if $g(X\cup Y)+g(X\cap Y)\leq g(X)+g(Y)$.
\begin{lemma}\label{lem:submod}
Let $g:2^E\longrightarrow\NN$ be nondecreasing and submodular, with $g(\varnothing)=0$, on a finite set $E$. For any integer $k\geq0$, the sublevel complex $D_k(g)=\{F\subseteq E:g(F)\leq k\}$ satisfies
\[
 \Ht_i(D_k(g);B)=0\qquad(i\geq k).
\]
For $k<0$ the sublevel complex is void.
\end{lemma}
\begin{proof}
Use induction on $|E|$ simultaneously for all integers $k$. The assertion is clear for $E=\varnothing$ under the augmented conventions. If $g(\{v\})=0$ for some $v$, submodularity and monotonicity imply $g(F\cup\{v\})=g(F)$ for every $F$; hence $D_k(g)$ is a cone for $k\geq0$. Otherwise fix $v\in E$ with $t=g(\{v\})\geq1$. The deletion of $v$ is the sublevel of the restricted function at level $k$. Its link is the sublevel at level $k-t$ of the function on $E\setminus\{v\}$ defined by
\[
 g_v(F)=g(F\cup\{v\})-t.
\]
This function is nonnegative, normalized, nondecreasing and submodular. There is a short exact sequence of augmented chain complexes for deletion and link,
\[
 0\to C_*(\del_vD_k;B)\to C_*(D_k;B)
 \to C_{*-1}(\lk_vD_k;B)\to0,
\]
with the shifted differential carrying the usual minus sign. This follows by separating the simplices containing $v$ from the others; it remains valid if $v$ is not a vertex of $D_k$. For $i\geq k$, the deletion homology in degree $i$ and the link homology in degree $i-1\geq k-t$ vanish by induction (or by voidness). The resulting long exact sequence proves the assertion.
\end{proof}
The bound is sharp; see Section~\ref{subsec:sharpness}.

\begin{lemma}[Local vanishing]\label{lem:local}
If $z$ is exact at level $c\geq1$, then
\begin{equation}\label{eq:local}
 \Ht_i(\G_{c-1}[z];B)=0\qquad(i\geq c-1).
\end{equation}
\end{lemma}
\begin{proof}
The function $L$ restricted to $2^z$ is submodular: all subsets of $z$ are independent, and $a(N(Y))$ is a weighted coverage function, whereas $a(Y)$ is modular. Set
\[
 g(F)=\min_{F\subseteq Y\subseteq z} L(Y).
\]
This function is integer-valued and nondecreasing. If $X,Y$ are minimizers for $F,F'$ respectively, their union and intersection are admissible candidates for $F\cup F'$ and $F\cap F'$. Hence
\[
 g(F\cup F')+g(F\cap F')\leq L(X\cup Y)+L(X\cap Y)\leq g(F)+g(F'),
\]
so $g$ is submodular. By Lemma~\ref{lem:internal}, the complex in question is $\{F:g(F)\leq c-1\}$. If there is a nonempty minimizer $M$ of $L$ on $2^z$, then
\[
 L(Y\cup M)\leq L(Y)+L(M)-L(Y\cap M)\leq L(Y).
\]
Thus each local witness may be extended by $M$, so the complex is a cone on any vertex of $M$. In the remaining case, the empty set is the unique minimizer, and $g(\varnothing)=0$, $g(F)>0$ for $F\ne\varnothing$. Lemma~\ref{lem:submod} at level $c-1$ proves \eqref{eq:local}.
\end{proof}
Neither flagness of $\G_1$ nor simplex structure of $\G_{c-1}[z]$ is needed, and both can fail; see Section~\ref{subsec:sharpness}.

\section{The negative and zero levels}\label{sec:base}
For a finite poset $P$, write $\mathcal O(P)$ for its order complex. We use the upper-fiber form of Quillen's Theorem A \cite[Proposition 1.6]{Quillen}: if $f:P\to Q$ is monotone and every upper fiber $\{p\in P:q\leq f(p)\}$ has a nonempty contractible order complex, then the induced map $\mathcal O(P)\to\mathcal O(Q)$ is a homotopy equivalence. Pointwise comparable monotone maps induce homotopic maps on order complexes. In particular, a monotone map pointwise below both the identity and a constant map provides a contraction.

Let $P_c=\{Z\in\Ind(G):Z\ne\varnothing,\ L(Z)\leq c\}$, ordered by inclusion. Whenever $\G_c$ has vertices, the inclusion of $P_c$ in its poset of nonempty faces induces a homotopy equivalence. To see this, fix a nonempty face $F$ and a minimizer $M\supseteq F$ with $L(M)=\lambda(F)$. For $Y\in P_c$ containing $F$, also $Y\vee M\supseteq F$. Therefore
\[
 L(Y\cap M)\leq L(Y)+L(M)-L(Y\vee M)\leq L(Y)\leq c.
\]
The map $Y\mapsto Y\cap M$ on the upper fiber is monotone and lies below both the identity and the constant map to $M$. The fiber is thus contractible and the fiber theorem applies.

\begin{lemma}[Low levels]\label{lem:negative}
For any graph and nonnegative integer weights, if $c<0$ then $\G_c$ is either void or contractible with vertices. If $\mu_G<0$, then $\G_0$ is contractible.
\end{lemma}
\begin{proof}
If $c<0$ and $\G_c$ is nonvoid, the least global minimizer $K_0$ is nonempty and belongs to $P_c$. For $Z\in P_c$, \eqref{eq:bi} gives $L(Z\cap K_0)\leq L(Z)\leq c<0$. Thus intersection with $K_0$ contracts $P_c$ by the same pointwise comparison, proving the first assertion.

For $c=0$ and $\mu_G<0$, let $P'_0$ be the subposet of $P_0$ consisting of sets meeting $K_0$. It contracts by intersection with $K_0$. If $Z\in P_0\setminus P'_0$, inequality \eqref{eq:bi} yields $L(Z\vee K_0)\leq L(Z)+\mu_G\leq\mu_G$. Hence $Z\vee K_0$ is a global minimizer and contains $K_0$. No vertex of $K_0$ has been removed; therefore $Z$ and $K_0$ are anticomplete and $M=Z\cup K_0$ is an independent global minimizer. For every $Y\in P'_0$ containing $Z$, intersection with $M$ remains in that upper fiber, since it contains $Z$, meets $K_0$, and by \eqref{eq:bi} has defect at most $L(Y)$. The fiber contracts to $M$. For $Z\in P'_0$ the upper fiber has least element $Z$. Applying the fiber theorem to $P'_0\subseteq P_0$ shows that $\mathcal O(P_0)$ is contractible, and therefore so is $\G_0$.
\end{proof}
In particular, for every $c<0$, the pair $(\G_c,\G_{c-1})$ is acyclic in augmented relative homology, by the exact sequence of the pair.

\begin{lemma}[Zero-level classification]\label{lem:zero}
Suppose $G$ is gap-free and $A>0$. Then $\G_0$ has no vertices, is contractible, or is homotopy equivalent to $S^0$. The last possibility can occur only if $A$ is even.
\end{lemma}
\begin{proof}
Lemma~\ref{lem:negative} handles $\mu_G<0$, so assume $\mu_G=0$.

\emph{Step 1: positive weights and clones.} First let every $a_v>0$. Replace each vertex $v$ by $a_v$ pairwise nonadjacent clones, with complete adjacency between two clone classes precisely when the original vertices are adjacent. Give each clone unit weight. Gap-freeness is preserved: an induced pair of disjoint edges cannot have two of its vertices in the same clone class, since those vertices have identical neighborhoods, and four distinct classes would project to an induced pair of disjoint edges in $G$. Saturating an independent clone set by taking all clones of its projected vertices does not change its neighborhood and cannot increase its defect. Hence $\mu_G$ is preserved, and a clone set belongs to the new $\G_0$ if and only if its projection belongs to the original $\G_0$. On nonempty face posets, saturation is an extensive order retraction onto the saturated faces, which correspond to the original nonempty faces. The homotopy type and total mass are thus preserved, including the case of no vertices. We may work with $n=A$ unit-weight vertices.

For every $X\subseteq V$, $Z=X\setminus N(X)$ is independent and $N(Z)\subseteq N(X)\setminus X$. Thus
\[
 |N(X)|-|X|\geq |N(Z)|-|Z|\geq 0.
\]
Hall's theorem \cite{Hall} in the bipartite double cover of $G$ gives a permutation $\sigma$ of $V$ such that $v\sigma(v)\in E$ for all $v$. Its cycles have length at least two. Call an independent set \emph{tight} if its defect is zero; these sets generate $\G_0$.

\emph{Step 2: an odd cycle.} If $\sigma$ has an odd cycle, consider any tight set $Z$. The two sets $\sigma(Z)$ and $\sigma^{-1}(Z)$ are contained in $N(Z)$ and have cardinality $|Z|=|N(Z)|$. Therefore
\[
 \sigma(Z)=N(Z)=\sigma^{-1}(Z),\qquad \sigma(N(Z))=Z.
\]
Every cycle of $\sigma$ meeting $N[Z]$ alternates between $Z$ and $N(Z)$ and is even. The fixed odd cycle is consequently outside $N[Z]$ for every tight $Z$. An edge between two tight sets, together with an edge of that odd cycle, would induce a gap. Thus the union of two tight sets is independent. By submodularity on this independent union and the inequality $L\geq0$, their union is tight. The union of all tight sets is tight, so $\G_0$ is a simplex, possibly with no vertices. This in particular covers odd $n$.

\emph{Step 3: a perfect matching and implications.} It remains to treat a permutation having only even cycles. Alternating edges on these cycles gives a perfect matching (a two-cycle contributes its one undirected edge). Write $\pi$ for the involution pairing its endpoints. Form a directed graph on $V$ by the implications
\[
 u\longrightarrow \pi(w) \quad\text{whenever } uw\in E\text{ is not a matching edge},
\]
including both directions of the original undirected edge in this rule. This is the usual implication-digraph structure of a $2$-CNF constraint system (see Aspvall--Plass--Tarjan \cite{APT}); the additional topological argument below uses gap-freeness. Let $u\Rightarrow v$ mean reachability, including the trivial path. The involution reverses all implications, so
\begin{equation}\label{eq:dual}
 u\Rightarrow v \quad\Longleftrightarrow\quad \pi(v)\Rightarrow\pi(u).
\end{equation}
Call a set consistent if it contains no matched pair. Tight sets are exactly the consistent sets closed under implications: for a tight set, $N(Z)=\pi(Z)$ follows by cardinality, and this gives closure. Conversely, a consistent closed set contains no matching edge, and cannot contain an edge $uw$ outside the matching, since closure would force $\pi(w)$. It is independent and satisfies $N(Z)=\pi(Z)$.

A vertex is valid if its reachability set is consistent. The vertices of $\G_0$ are exactly the valid ones. A set of valid vertices is a face if and only if the union of their reachability sets is consistent. Thus pairwise compatibility suffices for a face: every inconsistency would occur within one or between two reachability sets. For valid $u,v$, incompatibility is equivalent to $u\Rightarrow\pi(v)$. In fact, if some $t$ is reached from $u$ and $\pi(t)$ from $v$, duality gives $t\Rightarrow\pi(v)$; the converse is immediate.

If distinct valid vertices $u,v$ satisfy $u\Rightarrow v$, then they are compatible. Every face containing $u$ can have $v$ adjoined: an incompatible $x$ for $v$ would satisfy $x\Rightarrow\pi(v)\Rightarrow\pi(u)$, contradicting compatibility with $u$. It follows that the simplicial map sending $u$ to $v$ and fixing all other vertices is a retraction to the deletion of $u$, contiguous to the identity. In particular, it preserves homotopy type.

The valid vertices are closed under reachability. In the directed graph on them, retain one representative of every terminal strongly connected component. Each other valid vertex reaches a retained representative, and may be deleted by the preceding retraction. Pairwise compatibility is unchanged in the induced subcomplex on the remaining vertices, so these deletions can be made in any order while retaining all representatives. Terminal components of the valid subgraph are also terminal in the full directed graph, because successors of valid vertices are valid.

Suppose retained representatives $r_1,r_2$ in distinct terminal components $U,W$ are incompatible. Then $\pi(r_2)\in U$, and duality implies $U=\pi(W)$: the involution maps strongly connected components bijectively to strongly connected components. Since both are terminal, duality also makes both initial. The union $U\cup W$ is closed under $\pi$ and has no directed arcs between it and the outside. It follows that no edge of $G$ joins it to the outside: matching edges stay within, and any other such edge would yield an arc leaving the union. The incompatibilities among retained representatives therefore form disjoint pairs of dual terminal components. There is at most one such pair, since two distinct pairs would give mutually anticomplete matching edges, an induced gap. As the reduced complex is determined by pairwise compatibility, it is a simplex, a join of a nonempty simplex with $S^0$, just $S^0$, or has no vertices. The simplex and the join with a nonempty simplex are contractible. This proves the classification with positive weights, and the $S^0$ alternative requires even mass.

\emph{Step 4: zero weights.} If zero weights occur and an edge joins the two weight supports, choose its zero-weight endpoint $w$; then $a_w=0$ and $r=a(N(w))>0$. For every integer $k$,
\begin{align}
 \del_w\G_k(G,a)&=\G_k(G-w,a|_{V\setminus\{w\}}),\label{eq:delzero}\\
 \lk_w\G_k(G,a)&=\G_{k-r}(G-N[w],a|_{V\setminus N[w]}).\label{eq:linkzero}
\end{align}
For the deletion identity, removing $w$ from any independent witness cannot increase its defect, and deleting a zero-weight vertex does not change the neighborhood mass of sets avoiding it. The link identity follows from \eqref{eq:link}. At level zero, if $w$ is a face, its link is nonvoid at a negative level and is contractible by Lemma~\ref{lem:negative}. Attaching its star to the deletion preserves homotopy type; if $w$ is not a face, the deletion is the entire complex. Induct on the number of zero-weight vertices; deletion preserves the mass and gap-freeness. The remaining case is precisely that both weight supports are nonempty and no edge joins the two parts. In this case, one part must be edgeless, as otherwise there is an induced gap. Its vertices are isolated in $G$. Any isolated vertex can be adjoined to every witness without increasing the defect. Since $\G_0$ contains the empty face, it is then a cone. This completes the induction.
\end{proof}

Recall that $\epsilon(A,c)=0$ if $A+c$ is odd, and $\epsilon(A,c)=1$ if $A+c$ is even. For $A>0$, Lemmas~\ref{lem:negative} and \ref{lem:zero} imply
\begin{equation}\label{eq:relativezero}
 H_j(\G_0,\G_{-1};B)=0\qquad(j\geq\epsilon(A,0)).
\end{equation}
Indeed, when $\mu_G<0$ both complexes are contractible (integrality gives $\mu_G\leq-1$); otherwise $\G_{-1}$ is void and the assertion follows from the classification, including $\{\varnothing\}$.

\section{Relative vanishing and the parity bound}\label{sec:relative}
\begin{theorem}[Relative vanishing]\label{thm:relative}
Let $G$ be gap-free, $a\in\NN^V$, $c\geq1$ an integer, and $A\geq c+1$. For every abelian group $B$,
\begin{equation}\label{eq:relative}
 H_j(\G_c,\G_{c-1};B)=0\qquad(j\geq c+\epsilon(A,c)).
\end{equation}
\end{theorem}
\begin{proof}
We first assume that every weight is positive. In this part, $c\geq1$ suffices; the mass assumption $A\geq c+1$ is not needed. Start with $\sd\G_{c-1}$ and insert the nonempty faces of $\G_c\setminus\G_{c-1}$ as vertices of the face poset order complex in decreasing cardinality, using any order for ties. At the insertion of $F$, its level satisfies $\lambda(F)=c$. Every strictly larger face of $\G_c$ containing $F$ is new (by monotonicity of $\lambda$) and has already been inserted. Among proper nonempty subfaces of $F$, precisely those in $\G_{c-1}$ are present. Thus the link of the inserted vertex is the join
\begin{equation}\label{eq:insert}
 D_F=\sd(\G_{c-1}[F])*\sd(\lk_F\G_c),
\end{equation}
where the upper factor is identified by removing $F$. By \eqref{eq:link}, it equals $\sd(\G_{c-L(F)}(H_F))$. Both factors contain the empty face; for the lower factor, $F$ itself is not in $\G_{c-1}$. Every lower element precedes every upper element, so the join accounts for all chains in the link. The insertion attaches a cone on $D_F$; the relative chain groups of that step identify with augmented chains of $D_F$ shifted by one. In particular its relative homology in degree $j\geq1$ is $\Ht_{j-1}(D_F;B)$.

If $L(F)>c$, then by \eqref{eq:lambdaext} the upper factor is a nonvoid negative-level complex (its level equals the residual minimum). Lemma~\ref{lem:negative} makes it contractible with vertices. Thus the join is contractible, even when the lower factor is $\{\varnothing\}$. All the required relative groups of this step vanish.

In the other case, $L(F)=\lambda(F)=c$. The lower factor has no reduced homology in degrees $i\geq c-1$ by Lemma~\ref{lem:local}. The upper factor is a zero-level complex, with minimum zero and mass
\[
 A_F=A-2a(F)-c\equiv A+c\pmod{2}.
\]
If $A_F=0$, positivity forces $H_F$ to have no vertices, so its zero-level complex is $\{\varnothing\}$. If $A_F>0$, Lemma~\ref{lem:zero} applies because $H_F$ is an induced gap-free graph. For odd $A+c$, the upper factor is contractible or $\{\varnothing\}$. For even $A+c$, it may also be homotopy equivalent to $S^0$. A join with a nonempty contractible factor is contractible; the join with $\{\varnothing\}$ leaves the other factor unchanged, and a join with a complex homotopy equivalent to $S^0$ is a suspension up to homotopy (when the lower factor has no vertices, it is just the upper factor). Consequently
\[
 \Ht_i(D_F;B)=0\qquad(i\geq c-1+\epsilon(A,c)).
\]
This also covers the case of two factors equal to $\{\varnothing\}$. Hence the relative homology of each insertion vanishes in the range of \eqref{eq:relative}. The long exact sequence of a triple, applied successively to the filtration, shows that the same is true for the pair $(\sd\G_c,\sd\G_{c-1})$. Barycentric subdivision preserves homology of pairs, proving the result for positive weights.

Now allow zero weights. We use lexicographic induction on $(c,|\{v:a_v=0\}|)$ for $c\geq1$. The assertion for positive weights has just been proved, independently of the inductive hypothesis. The relative groups in negative levels vanish by Lemma~\ref{lem:negative}, while at level zero with positive mass we use \eqref{eq:relativezero}. Suppose $w$ is a zero-weight vertex with $r=a(N(w))>0$. Put $D_k=\del_w\G_k$ and $E_k=\lk_w\G_k$. There is a short exact sequence of augmented chain complexes
\begin{equation}\label{eq:chainpair}
 \begin{split}
 0\longrightarrow C_*(D_c,D_{c-1};B)&\longrightarrow C_*(\G_c,\G_{c-1};B)\\
 &\longrightarrow C_{*-1}(E_c,E_{c-1};B)\longrightarrow0.
 \end{split}
\end{equation}
To verify exactness, separate faces of the relative chain basis according as they avoid or contain $w$. Removing $w$ from the latter gives precisely the relative basis of the link pair. Orient $w$ first; the differential term removing $w$ lies in the deletion subcomplex, and all other terms give the shifted differential with a minus sign. These identifications work over any $B$, even when a link is void or has only the empty face.

The deletion pair has the same mass and fewer zero weights by \eqref{eq:delzero}, hence its relative homology in degree $j\geq c+\epsilon(A,c)$ vanishes by induction. By \eqref{eq:linkzero}, the link pair has levels $c'=c-r$ and $c'-1$ in $G-N[w]$, of mass $A'=A-r$. This is where the mass hypothesis is used: it passes to the link pair and gives $A'>0$ when $c'=0$. Indeed,
\[
 A'-c'=A-c\geq1,\qquad A'+c'\equiv A+c\pmod2.
\]
If $c'\geq1$, then $c'<c$ and $j-1\geq c'+\epsilon(A',c')$, so the induction hypothesis annihilates its relative homology in degree $j-1$. If $c'=0$, then $A'>0$ and $j-1\geq\epsilon(A',0)$, so \eqref{eq:relativezero} applies. If $c'<0$, negative-level relative homology is zero in every degree. The long exact sequence of \eqref{eq:chainpair} proves \eqref{eq:relative} in each case.

Finally, if both weight supports are nonempty but no edge joins the zero-weight and positive-weight parts, gap-freeness implies an isolated vertex as in Lemma~\ref{lem:zero}. The mass condition ensures that the positive part is nonempty. Since $c-1\geq0$, both levels contain the empty face and are cones on that isolated vertex. Their relative homology vanishes. This closes the induction and proves the theorem.
\end{proof}
\begin{remark}\label{rem:masszero}
The preceding treatment of zero weights is essential. A residual graph of zero mass need not be empty, and its zero-level complex need not satisfy Lemma~\ref{lem:zero}. See Section~\ref{subsec:sharpness}.
\end{remark}

\begin{theorem}[Parity-sensitive absolute vanishing]\label{thm:C}
Let $G$ be gap-free, $a\in\NN^V$, and $c\in\ZZ$ with $A\geq c+1$. Then, for every abelian group $B$,
\begin{equation}\label{eq:C}
 \Ht_j(\G_c(G,a);B)=0\qquad(j\geq c+\epsilon(A,c)).
\end{equation}
\end{theorem}
\begin{proof}
For $c<0$ apply Lemma~\ref{lem:negative}; for $c=0$ the mass hypothesis permits Lemma~\ref{lem:zero}. For $c\geq1$ we induct on $c$; the two parities alternate. The mass hypothesis holds also for $c-1$, and $\epsilon(A,c-1)=1-\epsilon(A,c)$. When $\epsilon(A,c)=0$, the induction hypothesis at level $c-1$ gives vanishing in degrees $j\geq c$. When $\epsilon(A,c)=1$, it gives vanishing in degrees $j\geq c-1$, in particular for $j\geq c+1$. The exact sequence of the pair together with Theorem~\ref{thm:relative} transfers the vanishing to $\G_c$ in the required range.
\end{proof}
The gap-free hypothesis cannot be omitted (Section~\ref{subsec:sharpness}). The odd-parity assertion of Theorem~\ref{thm:C} is the bound used for the algebraic application. Its even-parity companion is needed in the induction. The relative result also gives an isomorphism on reduced homology induced by $\G_{c-1}\hookrightarrow\G_c$ for $j\geq c+\epsilon(A,c)$, and an injection in degree $c+\epsilon(A,c)-1$, under the assumptions of Theorem~\ref{thm:relative}.

\section{Capacitated matchings and degree complexes}\label{sec:algebra}
For a nonnegative integer vector $b$ on a graph $H$, let $\nu_b(H)$ be the largest size of a multiset of edges whose incidences at each vertex $v$ number at most $b_v$. Put $|b|=b(V(H))$ and use $\mu_H$ with respect to these weights. Directly from the definition of an ordinary ideal power,
\begin{equation}\label{eq:member}
 x^b\in I(H)^q\quad\Longleftrightarrow\quad \nu_b(H)\geq q.
\end{equation}
\begin{lemma}[Capacitated matching identity]\label{lem:matching}
If $H$ is gap-free, then
\begin{equation}\label{eq:match}
 \nu_b(H)=\left\lfloor\frac{|b|+\mu_H}{2}\right\rfloor.
\end{equation}
\end{lemma}
\begin{proof}
Removing zero-weight vertices preserves both $\nu_b$ and $\mu_H$: removal of such vertices from an independent set cannot increase its defect, and the defect of an independent set of the positive support is unchanged. Replace each remaining vertex by $b_v$ unit clones as in Lemma~\ref{lem:zero}. Saturation preserves the minimum defect, and distributing incidences among distinct clones identifies the maximum capacitated matching size with the ordinary matching number of the clone graph. This graph is gap-free. We may therefore assume unit weights on $n=|b|$ vertices, with the empty case immediate.

Write $\nu$ for the matching number and $\delta=n-2\nu$ for its deficiency. An independent set $Z$ can have at most $|N(Z)|$ of its vertices matched, so $\delta\geq -\mu_H$. By the Tutte--Berge formula (see \cite{LP} or \cite[Section 24.1]{Schrijver}), there is $D\subseteq V(H)$ with $\delta=\operatorname{odd}(H-D)-|D|$, where $\operatorname{odd}$ counts odd-order components. Gap-freeness implies that $H-D$ has at most one component containing an edge. Let $Z$ consist of its isolated vertices. Then $N_H(Z)\subseteq D$, whence
\[
 \delta\leq |Z|+1-|D|\leq |Z|-|N_H(Z)|+1\leq -\mu_H+1.
\]
Since $\delta\equiv n\pmod2$, these two bounds determine $\delta=n-2\lfloor(n+\mu_H)/2\rfloor$, proving \eqref{eq:match}.
\end{proof}
Lemma~\ref{lem:matching} is a gap-free form of the Tutte--Berge formula. It includes the edgeless case; see Section~\ref{subsec:sharpness} for a counterexample without gap-freeness.

\paragraph{Fractional interpretation.}
Let $\nu_b^*(H)$ denote the maximum size of a fractional $b$-matching. Its objective is $\sum_{e\in E(H)}y_e$, subject to $y_e\geq0$ and $\sum_{e\ni v} y_e\leq b_v$. For arbitrary graphs,
\begin{equation}\label{eq:fractional}
 2\nu_b^*(H)=|b|+\mu_H;
 \qquad \nu_b(H)=\lfloor\nu_b^*(H)\rfloor\quad\text{if $H$ is gap-free}.
\end{equation}
The first identity follows from a maximum-flow computation on the bipartite double cover. Take left and right copies of $V(H)$, arcs $u_L\to v_R$ for each directed incidence of an edge, and source and sink arcs of capacities $b_u,b_v$, respectively. Give middle arcs capacity greater than $|b|$. The maximum flow equals $|b|+\min_{X\subseteq V}(b(N(X))-b(X))$ by the minimum cut formula: for any left subset $X$, the right source-side vertices must include $N(X)$ in a minimum cut. The minimum over all $X$ is $\mu_H$. Indeed $Z=X\setminus N(X)$ is independent, $N(Z)\subseteq N(X)\setminus X$, and $b(N(X))-b(X)\geq b(N(Z))-b(Z)$. The reverse inequality for minima follows by taking independent $X$. Every fractional matching gives a feasible flow of twice its size by assigning $y_{uv}$ in each direction. Conversely, averaging the two directed middle-arc flows on each edge yields a feasible fractional matching of half the flow value; its vertex incidence sum is at most $b_v$. Thus the maximum flow is exactly $2\nu_b^*(H)$, including the case $|b|=0$. The second identity follows from Lemma~\ref{lem:matching}. This is the fractional $b$-matching relaxation without individual upper bounds on edge variables (compare \cite[Section 30.2 and Chapters 21, 31]{Schrijver}). This exact rounding in the gap-free case is what makes the levels in the degree-complex formula integral.

For a monomial ideal $J\subset S$ and $\gamma\in\ZZ^V$, let $\eta=\{v:\gamma_v<0\}$, and write $\gamma^+_v=\max(\gamma_v,0)$. The multigraded degree complex, in the localization formulation of Minh--Trung \cite[Lemma 1.2]{MT}, is
\begin{equation}\label{eq:degree-def}
 \Delta_\gamma(J)=\{F\subseteq V\setminus\eta:
 x^\gamma\notin J S[x_v^{-1}:v\in F\cup\eta]\}.
\end{equation}
We recall the localization criterion for powers in a form that also handles negative coordinates. For an independent set $W$, set $H_W=G-N_G[W]$. In the ring obtained by inverting the variables of $W$,
\[
 I(G)S[x_W^{-1}]=\bigl((x_v:v\in N_G(W))+I(H_W)\bigr)S[x_W^{-1}].
\]
Indeed, each edge incident to $W$ yields a neighbor variable; other edges meeting a neighbor are redundant; the remaining edges are those of $H_W$. The variable sets of the two summands on the right are disjoint. Localization commutes with powers. For $W\supseteq\eta$, the expansion of the $q$th power of the sum and \eqref{eq:member} therefore give
\begin{equation}\label{eq:criterion}
 \begin{gathered}
 x^\gamma\in I(G)^q S[x_W^{-1}]\\
 \Longleftrightarrow\quad
 \gamma^+(N_G(W))+\nu_{\gamma^+|_{V(H_W)}}(H_W)\geq q.
 \end{gathered}
\end{equation}
For the first summand, a product of neighbor variables can have length at most their available total exponent; for the second, the available length is exactly the maximum capacitated matching. If $W$ is not independent, an edge becomes a unit and every monomial belongs to the localized ideal. This criterion is the localization form of the neighborhood and monomial-order test used in \cite[Lemma 2.22]{MV}; we include its justification to track all multigrades.

\begin{proposition}[Degree-complex identity]\label{prop:degree}
Let $G$ be gap-free and $q\geq1$. For $a\in\NN^V$ with $A=a(V)$,
\begin{equation}\label{eq:degree}
 \Delta_a(I(G)^q)=\G_{2q-1-A}(G,a).
\end{equation}
More generally, for $\gamma\in\ZZ^V$ and negative support $\eta$, the complex is void if $\eta$ is not independent. If $\eta$ is independent, put
\[
 H=G-N_G[\eta],\qquad b=\gamma^+|_{V(H)},\qquad
 r=\gamma^+(N_G(\eta)),\qquad q'=q-r.
\]
If $q'\leq0$ the complex is void. Otherwise,
\begin{equation}\label{eq:negative-degree}
 \Delta_\gamma(I(G)^q)=\Delta_b(I(H)^{q'})=\G_{2q'-1-|b|}(H,b),
\end{equation}
where, as above, $|b|=b(V(H))$.
\end{proposition}
\begin{proof}
For $a\geq0$ and $F$ independent, the criterion \eqref{eq:criterion}, followed by Lemma~\ref{lem:matching}, yields nonmembership in the localization exactly when
\[
 2a(N(F))+a(V(H_F))+\mu_{H_F}\leq2q-1.
\]
By \eqref{eq:lambdaext} and \eqref{eq:mass}, the left side is $A+\lambda_G(F)$. Nonindependent sets are not faces. This proves \eqref{eq:degree}, including the empty face and void cases.

For a general $\gamma$, first invert the variables of $\eta$. If it contains an edge, the ideal becomes the whole ring. Otherwise the neighbor variables generate the first summand as above. When $r\geq q$, the monomial already belongs to the localized power. If $r<q$, further inverting any variable in $N_G(\eta)$ makes the ideal the whole ring, so any additional face must be contained in $V(H)$. Localizing further at $F\subseteq V(H)$ preserves the first neighbor-variable summand and its available mass $r$. Expanding the power in these disjoint variable sets shows that nonmembership is exactly nonmembership of $x^b$ in the localized $I(H)^{q-r}$. This gives the first equality of \eqref{eq:negative-degree}, and \eqref{eq:degree} gives the second, since $H$ is gap-free. All powers of residual edge ideals here have positive exponents, including when the residual ideal is zero.
\end{proof}

\section{Proof of the main theorem}\label{sec:proof}
For a nonzero proper monomial ideal $J$, Takayama's formula \cite{Tak} in the localization formulation \cite[Section 1]{MT} expresses the multigraded local cohomology of $S/J$ in terms of \eqref{eq:degree-def}. For fixed negative support $\eta$ the complex is independent of the actual negative entries:
\begin{equation}\label{eq:tak}
 \dim_K H^i_{\m}(S/J)_\gamma=
 \dim_K\Ht_{i-|\eta|-1}(\Delta_\gamma(J);K).
\end{equation}
Here $\Delta_\gamma$ is used with the augmented conventions of Section~\ref{sec:prelim}; the formula follows from the $\gamma$-graded \v{C}ech complex, whose nonzero summands correspond to inverted sets $W=F\cup\eta$ with $F\in\Delta_\gamma(J)$. Their cochain degree is shifted by $|\eta|+1$. Over a field, reduced cohomology and homology have equal dimension. This explains the contributions of $\{\varnothing\}$ and the void complex as well.

The local cohomological description of regularity gives
\begin{equation}\label{eq:reglocal}
 \reg J=1+\max_{\substack{i\geq0,\ \gamma\in\ZZ^V\\ H^i_{\m}(S/J)_\gamma\ne0}}\left(\sum_v\gamma_v+i\right).
\end{equation}
With the negative support fixed, the degree complex is unchanged by replacing each negative coordinate by $-1$, and this maximizes the coordinate sum. Hence \eqref{eq:tak} gives
\begin{equation}\label{eq:reg}
 \reg J=\max\bigl\{|\gamma^+|+j+2:
 \Ht_j(\Delta_\gamma(J);K)\ne0\bigr\},
\end{equation}
where $\gamma_v\in\{-1\}\cup\NN$ for all $v$, and $|\gamma^+|=\sum_v\gamma_v^+$. Compare also \cite[Lemmas 2.9 and 2.11]{MV}. The maximum is attained, by the usual finiteness of Castelnuovo--Mumford regularity.

\begin{proof}[Proof of Theorem~\ref{thm:main}]
Write $R=\reg_K I(G)\geq2$. Consider a multidegree and an integer $j$ contributing to \eqref{eq:reg} for $J=I(G)^q$. By Proposition~\ref{prop:degree}, the negative support is independent, $q'=q-r\geq1$, and, setting $c'=2q'-1-|b|$,
\begin{equation}\label{eq:D}
 \Ht_j(\G_{c'}(H,b);K)\ne0,
 \qquad D=|\gamma^+|+j+2=r+|b|+j+2.
\end{equation}
We distinguish two exhaustive possibilities.

\emph{Low mass: $|b|\leq q'-1$.} Every independent set $Z$ in $H$ satisfies $L_H(Z)\leq|b|\leq c'$, so $\G_{c'}(H,b)=\Ind(H)$. Consider the first power $I(G)$ and the multidegree with entries $-1$ on the same negative support $\eta$ and zero elsewhere. By Proposition~\ref{prop:degree} its degree complex is precisely $\Ind(H)$. Using \eqref{eq:reg} at this degree yields $j+2\leq R$. Consequently
\begin{equation}\label{eq:low}
 D\leq r+q'-1+R=q+R-1.
\end{equation}
This use of low mass is consistent with the degree-complex observation in \cite[Lemma 2.15]{MV}.

\emph{High mass: $|b|\geq q'$.} In this case,
\[
 |b|+c'=2q'-1\text{ is odd},\qquad |b|-c'=2|b|-2q'+1\geq1.
\]
Theorem~\ref{thm:C} therefore applies to the gap-free residual $H$, with mass $|b|$, even when $c'$ is negative or some weights vanish. The nonvanishing in \eqref{eq:D} forces $j\leq c'-1$. Hence
\begin{equation}\label{eq:high}
 D\leq r+|b|+c'+1=r+2q'=2q-r\leq2q.
\end{equation}
Together \eqref{eq:low} and \eqref{eq:high} prove
\[
 \reg I(G)^q\leq\max\{2q,q+R-1\}.
\]

For the reverse inequality, the ideal is nonzero and generated in degree $2q$, so its regularity is at least $2q$. If $R=2$ this is sufficient. Suppose $R>2$. Choose a degree attaining the maximum in \eqref{eq:reg} for $q=1$. Inequality \eqref{eq:high} shows that its contribution cannot arise from the high-mass case. The low-mass case at $q=1$ necessarily has $r=0$, $q'=1$ and $|b|=0$. Thus for some independent $\eta$,
\[
 \Ht_{R-2}(\Ind(G-N[\eta]);K)\ne0.
\]
The residual has a vertex $v$, because $R-2>0$. For arbitrary $q\geq1$, choose the degree with entries $-1$ on $\eta$, entry $q-1$ at $v$, and zero elsewhere. It has $r=0$ and $|b|=q-1$, so the same independence complex is obtained by the low-mass argument. Its nonzero homology contributes $q+R-1$ in \eqref{eq:reg}. Thus $\reg I(G)^q\geq q+R-1$, also agreeing with the general lower bound of \cite[Lemma 2.16]{MV}. The two bounds prove \eqref{eq:main}. As $I(G)^q$ is generated in degree $2q$, its resolution is linear exactly when its regularity is $2q$, equivalently $q\geq R-1$.
\end{proof}

\section{Consequences and scope}\label{sec:consequences}
For completeness, we include a short argument for the induced-gap obstruction of Francisco--H\`a--Van Tuyl \cite[Proposition 1.8]{NP}; the multidegree used below differs slightly from the one in \cite{NP}. If $G$ contains an induced pair of disjoint edges with corresponding monomials $m_1,m_2$, the degree $m_1^qm_2$ is divisible by exactly two minimal generators of $I(G)^q$, namely $m_1^q$ and $m_1^{q-1}m_2$. No generator involving another vertex can divide this degree. The two-generator syzygy in this multidegree is minimal: no proper divisor of the multidegree is divisible by both generators. Its degree is $2q+2$, so $\reg I(G)^q\geq2q+1$.

\begin{corollary}\label{cor:np}
For a finite simple graph $G$ with at least one edge, the following are equivalent over any field $K$:
\begin{enumerate}[label=\textup{(\roman*)}]
 \item $G$ is gap-free;
 \item $I(G)^q$ has a linear free resolution for some integer $q\geq1$;
 \item $I(G)^q$ has a linear free resolution for all sufficiently large $q$.
\end{enumerate}
Under these equivalent conditions, the least such exponent is $\reg_K I(G)-1$, and for every $q\geq1$,
\[
 \reg_K I(G)^{q+1}=\max\{2q+2,\reg_K I(G)^q+1\}.
\]
\end{corollary}
\begin{proof}
Theorem~\ref{thm:main} gives (i)$\Rightarrow$(iii)$\Rightarrow$(ii) and the exact threshold. The induced-gap obstruction gives (ii)$\Rightarrow$(i). For the recurrence, substitute \eqref{eq:main} and note that $2q+2>2q+1$. In particular the inequality proposed in \cite[Open Problem 1.11(2)]{NP} holds with equality.
\end{proof}
\begin{corollary}\label{cor:linear-power}
For an integer $q\geq1$, the ideal $I(G)^q$ has a linear free resolution if and only if $G$ is gap-free and $\reg_K I(G)\leq q+1$.
\end{corollary}
\begin{proof}
Use the induced-gap obstruction and Theorem~\ref{thm:main}.
\end{proof}
For $q=2$ and $q=3$, Corollary~\ref{cor:linear-power} is the theorem of Minh and Vu \cite{MV}.
\subsection{Upper bounds and the stabilization index}
For an equigenerated ideal of degree two let $b(I)$ denote the eventual intercept and $c(I)$ the least positive integer such that $\reg I^{q+1}-\reg I^q=2$ for all $q\geq c(I)$, following \cite[Section 4]{BN}. For edge ideals, $c(I(G))$ is the regularity stabilization index $\operatorname{rstab}(G)$ of Minh and Vu \cite[p.~338]{MV}.
\begin{corollary}\label{cor:further}
Let $G$ be gap-free with an edge and write $R=\reg_K I(G)$. Then:
\begin{enumerate}[label=\textup{(\roman*)}]
 \item $\reg_K I(G)^q\leq2q+R-2$ for every $q\geq1$, with equality if and only if $q=1$ or $R=2$.
 \item If $R\leq3$, every ordinary power with exponent $q\geq2$ has a linear free resolution.
 \item $b(I(G))=0$, $c(I(G))=R-1$, and for $q\geq1$,
 \[
  \reg_K I(G)^{q+1}-\reg_K I(G)^q=
  \begin{cases}1,&1\leq q\leq R-2,\\2,&q\geq R-1.\end{cases}
 \]
 The first range is empty when $R=2$.
 \item $I(G)^2$ has a linear free resolution if and only if $R\leq3$.
\end{enumerate}
\end{corollary}
\begin{proof}
The identity
\[
 \max\{2q,R+q-1\}=2q+R-2-\min\{q-1,R-2\}
\]
proves the first assertion, including its equality cases. The second and fourth follow by comparing $R+q-1$ with $2q$. The third follows from the two ranges of \eqref{eq:main}, which meet at $q=R-1$; the increments before that exponent equal one and those at or after it equal two. This proves minimality in the definition of $c(I(G))$.
\end{proof}
Part (i) establishes the gap-free case of the Alilooee--Banerjee--Beyarslan--H\`a--Jayanthan--Selvaraja conjecture as formulated in \cite[Conjecture 1.1]{BN} (Conjecture A of \cite{BBH}). Its assertion for arbitrary graphs is not established here. Following \cite{BBH}, $G$ is locally of regularity at most $r-1$ if $\reg(I(G):x)\leq r-1$ for every vertex $x$; this implies $R\leq r$, as noted in \cite[Section 1]{BBH} (it also follows from \cite[Lemma 3.1]{DHS} by induction on the number of vertices). Thus part (i) also yields the gap-free bound of \cite[Theorem 4.2]{BBH}, and part (ii) recovers \cite[Theorem 4.5]{BBH}, which treats locally linear gap-free graphs. Part (ii) settles the case $\reg I(G)=3$, which Nevo and Peeva single out as open at the end of their discussion of \cite[Open Problems 1.11]{NP}, and answers \cite[Question 1.1]{Ban}, \cite[Conjecture B]{BBH}, and \cite[Conjecture 4.2]{BN} in their full range of exponents. The previously known cases $q=2,3$ are due to Minh and Vu \cite{MV}.

\begin{corollary}\label{cor:unbounded}
Over every fixed field, the stabilization indices $c(I(G))$ of gap-free edge ideals are unbounded. On an $n$-vertex gap-free graph, $c(I(G))=O(\log n)$, with an absolute implicit constant.
\end{corollary}
\begin{proof}
Januszkiewicz and \'{S}wi\k{a}tkowski \cite[Section 6]{JS} constructed flag orientable closed pseudomanifolds without empty squares in every dimension; see \cite{CKV} and \cite{BC} for their use in regularity questions. Let $\Delta$ be such a pseudomanifold of dimension $d\geq1$, and let $G$ be the complement of its one-skeleton. Since $\Delta$ is flag, $\Ind(G)=\Delta$; since $\Delta$ has no empty square, $G$ is gap-free; and $G$ has an edge because $\Delta$ is not a simplex. As $\Delta$ is orientable, $\Ht_d(\Delta;K)\ne0$ for every field $K$, and no induced subcomplex has homology above degree $d$. Hochster's formula (Section~\ref{sec:examples}) therefore gives $R=d+2$, and Corollary~\ref{cor:further} gives $c(I(G))=d+1$. For the upper bound, gap-free graphs are exactly those whose edge ideals have linear first syzygies \cite[Proposition 1.3]{NP}, \cite[Corollary 2.9]{DHS}, and the logarithmic bound of Dao, Huneke and Schweig \cite[Theorem 4.9]{DHS}, with $k=1$, gives $R=O(\log n)$ with absolute constants. Hence $c(I(G))=R-1=O(\log n)$.
\end{proof}
This gives an explicit answer for gap-free edge ideals to the stabilization-index question in \cite[Question 4.1]{BN}. Unboundedness here means that the index cannot be bounded independently of the number of vertices even with generator degree fixed at two. Minh and Vu \cite[pp.~338--339]{MV} recalled both ingredients, the logarithmic bound of \cite{DHS} and the gap-free graphs of arbitrary regularity of \cite{CKV}, and observed that their Conjecture A would allow arbitrarily large stabilization indices; Corollary~\ref{cor:unbounded} confirms this and determines the index exactly.

\paragraph{Relation to known cases.}
By Fr\"oberg \cite{Fro}, $R=2$ is equivalent to chordality of the complement, which implies gap-freeness. In that case Theorem~\ref{thm:main} recovers the edge-ideal case of Herzog--Hibi--Zheng \cite{HHZ}. Combined with the bounds $R\leq3$ of Nevo \cite{Nevo} (gap-free claw-free) and Banerjee \cite{Ban} (gap-free cricket-free), it gives linearity for every $q\geq2$ in those classes; it does not itself prove those initial-regularity bounds. Conversely, by the square criterion of Minh and Vu \cite{MV} (Corollary~\ref{cor:further}(iv)), the known linearity of the square for gap-free diamond-free graphs \cite{Erey18}, and for the families in \cite{E26,E3,D4} with linear powers from the second one onward, forces $R\leq3$ in those classes. The new uniform implication is that every gap-free graph with $R\leq3$ has linear powers for all exponents $q\geq2$. Theorem~\ref{thm:main} also yields the recurrence with equality in Corollary~\ref{cor:np}, for arbitrary initial regularity. All of these comparisons are consequences or consistency checks, not assumptions of the proof.

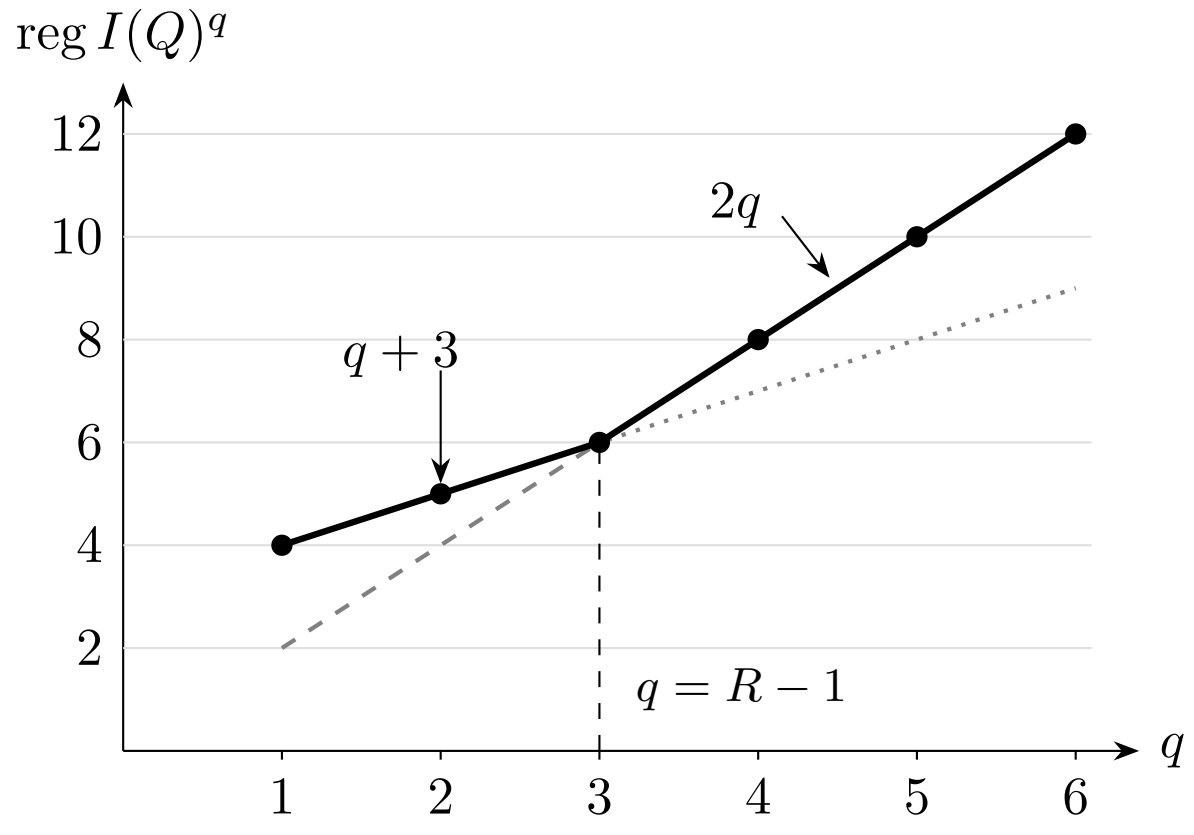
\begin{figure}[htbp]
\centering
\begin{tikzpicture}[x=1.05cm,y=.34cm,>=Stealth,font=\small]
\draw[->] (0,0)--(6.4,0) node[right] {$q$};
\draw[->] (0,0)--(0,13) node[above] {$\reg I(Q)^q$};
\foreach \x in {1,2,3,4,5,6}{\draw (\x,0)--(\x,-.17) node[below] {\x};}
\foreach \y in {2,4,6,8,10,12}{\draw[gray!25] (0,\y)--(6.1,\y);\node[left] at (0,\y){\y};}
\draw[gray, dashed, thick] (1,2)--(6,12);
\draw[gray, dotted, thick] (1,4)--(6,9);
\draw[black,very thick] (1,4)--(3,6)--(6,12);
\foreach \x/\y in {1/4,2/5,3/6,4/8,5/10,6/12}{\fill (\x,\y) circle (2pt);}
\draw[dashed] (3,0)--(3,6);\node[anchor=west,font=\footnotesize] at (3.1,1.2) {$q=R-1$};
\node[anchor=west] at (1.25,7.7) {$q+3$};\draw[thin,->] (2.0,7.4)--(2,5.2);
\node[anchor=east] at (4.15,10.6) {$2q$};\draw[thin,->] (4.15,10.4)--(4.45,9.2);
\end{tikzpicture}
\caption{Regularity of the powers of $I(Q)$, where $Q$ is the complement of the icosahedron ($R=4$). The dashed line is $2q$, the dotted line is $q+R-1=q+3$, and the solid line is their maximum; the dots are the values $\reg I(Q)^q$.}\label{fig:regcurve}
\end{figure}

\paragraph{Dependence on the field.}
Theorem~\ref{thm:C} holds with arbitrary abelian coefficients. Thus for ordinary powers all field dependence in \eqref{eq:main} is carried by $R=\reg_K I(G)$; the formula does not assert that $R$ is characteristic-independent. For a fixed $R$, the entire regularity function is the maximum of two affine functions, with the only possible change of slope at $q=R-1$ (Figure~\ref{fig:regcurve}).

\section{Examples and limitations}\label{sec:examples}
We use Hochster's formula (see \cite[Corollary 5.12]{MS}) in the form
\[
 \reg_K I(G)=\max\{j+2:\Ht_j(\Ind(G[W]);K)\ne0,\ W\subseteq V\}
\]
for graphs having an edge. The following examples illustrate the main theorem; none of them requires a separate computation of powers.

\paragraph{The five-cycle.}
For $G=C_5$ with unit weights, $A=5$. Every nonempty independent set has defect one: a singleton has two neighbors, and an independent pair has three neighbors. Hence $\mu=0$, $\Gamma_0=\{\varnothing\}$, and $\Gamma_1=\Ind(C_5)\simeq S^1$. The graph is gap-free and Hochster's formula gives $R=3$. Thus $\reg I^q=2q$ for all $q\geq2$, in agreement with the cycle formulas of \cite{BHT}. The unit-weight levels also show that the parity bound is attained: for $a=(1,\ldots,1)$, Proposition~\ref{prop:degree} gives $\Delta_a(I^3)=\Gamma_0=\{\varnothing\}$ and $\Delta_a(I^4)=\Gamma_2=\Ind(C_5)$, which contribute $5-1+2=6$ and $5+1+2=8$ to \eqref{eq:reg}. These are the values $2q$ allowed by \eqref{eq:high} for $q=3,4$.

\paragraph{The complement of the icosahedron.}
Take $Q$ to be the complement of the one-skeleton of the icosahedron, the graph in Nevo--Peeva's Counterexample 1.10 \cite{NP}. Its independence complex is the flag triangulation of $S^2$ given by the boundary of the icosahedron, and its complement has no induced square. Thus $Q$ is gap-free and $R=4$ over every field by Hochster's formula. Theorem~\ref{thm:main} gives the full sequence
\[
 (\reg I(Q)^q)_{q=1}^5=(4,5,6,8,10),\qquad \reg I(Q)^q=2q\quad(q\geq3).
\]
The values for $q=2,3,4$ agree with the Macaulay2 computations reported in \cite[Counterexample 1.10]{NP}, the last two due to M.~Stillman. The square is nonlinear and the cube is linear. The values for later powers are deduced from the theorem rather than presented as new numerical computations.

\paragraph{The dunce-hat construction.}
Banerjee and Nevo \cite[Section 4]{BN} describe the complement of the one-skeleton of a flag-no-square triangulation of the dunce hat. Their construction has $R=3$ and $\reg(I(G)^2:x_ux_v)=3$ for each edge $uv$. It shows why their direct colon-ideal estimate cannot by itself prove linearity of the square. Theorem~\ref{thm:main} gives $\reg I(G)^2=4$ and all higher powers linear. The square value also follows from the theorem of Minh--Vu \cite{MV}, so it is not claimed as a new isolated consequence. The point here is the uniform proof for arbitrary exponents, not an independent computation of a particular triangulation.

\subsection{Sharpness of the hypotheses}\label{subsec:sharpness}

\paragraph{Exactness and adjacent levels.}
Both hypotheses of Lemma~\ref{lem:internal} are needed. For exactness, take an edge $uv$ and an isolated vertex $w$ with respective weights $1,2,1$. Then $L(\{u\})=1$ but $\lambda(\{u\})=0$, so $\{u\}$ is not exact. Although $u\in\G_0$, no subset of $\{u\}$ containing $u$ has defect at most zero. For the adjacent level, take the graph on $0,1,2,3$ with edges $01,03,13,23$ and weights $(1,3,1,1)$. The set $z=\{0,2\}$ is exact at level two and $\{2\}\in\G_{-2}$, since $\lambda(\{2\})=-2$; but the minimum of $L(Y)$ over $\{2\}\subseteq Y\subseteq z$ is zero. Thus \eqref{eq:internal} fails if $c-1$ is replaced by a lower level.

\paragraph{Sharp submodular bound.}
The range in Lemma~\ref{lem:submod} is sharp: for $g(F)=|F|$ on $k+1$ elements, $D_k(g)$ is the boundary of a $k$-simplex, with reduced homology $B$ in degree $k-1$ (also for $k=0$ under our convention).

\paragraph{Nonflagness and nonsimplex local complexes.}
Flagness of $\G_1$ and the simplex property of $\G_{c-1}[z]$ can both fail for gap-free graphs, which is why Lemma~\ref{lem:local} is formulated homologically. Take an independent set $U=\{u_1,\ldots,u_m\}$ and a clique $P=\{p_1,\ldots,p_m\}$ with the only cross edges $u_ip_i$. Give each $u_i$ weight one and each $p_i$ weight two. This split graph is gap-free. For $Y\subseteq U$, $L(Y)=|Y|$. With $m=2$, $z=U$ is exact at level two (there is no further independent extension), and $\G_1[U]$ consists of two isolated vertices. For a nonflag example, take $m=3$ and adjoin an isolated vertex of weight one. Each pair of vertices in $U$ extends by the isolated vertex to an independent set of defect one, whereas $U$ itself can only be extended by that isolated vertex and has minimum defect two. In particular the three vertices of $U$ are pairwise joined by edges of $\G_1$ but do not span a face of $\G_1$.

\paragraph{Zero-mass residual.}
The separate treatment of zero weights in Theorem~\ref{thm:relative} is essential. A residual graph with zero total mass need not be empty. For instance, take the graph join of $K_2$ with the disjoint union of a vertex $z$ and $C_5$, giving weight one to the vertices of $K_2$ and to $z$, and zero to the pentagon. This graph is gap-free, $A=3$, and $L(z)=\lambda(z)=1$. The residual of $z$ is the zero-weight pentagon, whose zero-level complex is $\Ind(C_5)\simeq S^1$. Thus Lemma~\ref{lem:zero} cannot be applied to that residual. Sequence \eqref{eq:chainpair} avoids such an application.

\paragraph{Necessity of gap-freeness in the topological theorem.}
The gap-free hypothesis in Theorem~\ref{thm:C} cannot be omitted: for $G=2K_2$ with unit weights, every independent set has defect zero. Hence $\G_1=\Ind(G)=S^0*S^0\simeq S^1$, whereas $A=4$, $c=1$, and $\epsilon(A,c)=0$. The required vanishing in degree one fails over any nonzero coefficient group.

\paragraph{Integer rounding without gap-freeness.}
Two disjoint triangles with unit weights have $|b|=6$, $\mu_H=0$, and $\nu_b(H)=2$. Their fractional matching number is three, obtained by assigning $1/2$ to every edge. Thus the gap-free hypothesis in the second part of \eqref{eq:fractional} is essential.

\begin{center}
\small
\begin{tabularx}{\linewidth}{@{}p{.36\linewidth}X@{}}
\toprule
Example & Limitation illustrated\\
\midrule
Edge and isolated vertex; weighted four-vertex graph & Exactness and the adjacent level in Lemma~\ref{lem:internal}.\\
Cardinality sublevel & Sharp range of Lemma~\ref{lem:submod}.\\
Split graphs & Local complexes need not be simplices; $\Gamma_1$ need not be flag.\\
$K_2$ joined to $\{z\}\sqcup C_5$ & A zero-mass residual can have nontrivial homology.\\
$2K_2$ & Gap-freeness in Theorem~\ref{thm:C}.\\
$2K_3$ & Gap-freeness in Lemma~\ref{lem:matching}.\\
\bottomrule
\end{tabularx}
\end{center}

\section{Further questions}\label{sec:open}
We close with problems that lie beyond the methods of this paper.
\begin{enumerate}
\item \emph{Symbolic powers.} Establish the symbolic-power equality in \cite[Conjecture A]{MV} for general gap-free graphs.
\item \emph{Linear quotients.} Does $I(G)^q$ have linear quotients whenever $G$ is gap-free and $q\geq R-1$? This is a stronger requirement than linearity, related to the persistence question in \cite{E26}.
\item \emph{Field dependence.} Determine which gap-free graphs have initial regularity depending on the field, and thus which can exhibit characteristic dependence in the initial part of the regularity function. Torsion in induced flag-no-square subcomplexes is relevant, but must be compared with the maximal homology degrees across all induced subcomplexes.
\item \emph{Beyond gap-freeness.} Find formulas for the intercept and the stabilization index of general edge ideals, and prove the ABBHJS bound beyond this class. The integer rounding and parity arguments above require additional information when there is an induced gap.
\item \emph{Other ideals and topology.} Investigate analogues of the degree-complex model for hypergraph edge ideals, path ideals, and suitable quadratic ideals, and stronger structural properties of the Hall filtration. For the latter, any extension of parity-sensitive vanishing or discrete Morse and shellability approaches must account for the counterexamples in Section~\ref{subsec:sharpness}.
\end{enumerate}

\section*{Use of generative AI}
Two generative AI systems were used in this work: OpenAI's ChatGPT and Anthropic's Claude. They were used to explore proof strategies, to write and run computational verification scripts, to produce independent referee-style evaluations of successive drafts, and to draft and revise the manuscript. The author is responsible for all mathematical claims, references, and other contents of the paper.

\end{document}